\documentclass[12pt,twoside]{article}

\usepackage{amsfonts}
\usepackage{amssymb}
\usepackage{amsmath}
\usepackage{graphicx}

\newcommand{\ignore}[1]{}

\def\@begintheorem#1#2{\par\bgroup{\sc #1\ #2. }\it\ignorespaces}
\def\@opargbegintheorem#1#2#3{\par\bgroup{\sc #1\ #2\ (#3). } \it\ignorespaces}
\def\@endtheorem{\egroup}
\newtheorem{theorem}{Theorem}[section]
\newtheorem{corollary}[theorem]{Corollary}
\newtheorem{lemma}[theorem]{Lemma}
\newtheorem{example}[theorem]{Example}
\newtheorem{proposition}[theorem]{Proposition}
\newtheorem{question}[theorem]{Question}
\newtheorem{definition}[theorem]{Definition}
\newcommand{\bt}[1]{\begin{theorem}\label{#1}}
\newcommand{\bc}[1]{\begin{corollary}\label{#1}}
\newcommand{\bl}[1]{\begin{lemma}\label{#1}}
\newcommand{\be}[1]{\begin{example}\label{#1}}
\newcommand{\bp}[1]{\begin{proposition}\label{#1}}
\newcommand{\bq}[1]{\begin{question}\label{#1}}
\newcommand{\ba}[1]{\begin{algorithm}\rm\label{#1}}
\newcommand{\bd}[1]{\begin{definition}\rm\label{#1}}
\newcommand{\bpr}{\noindent {\em Proof. }}
\newcommand{\et}{\end{theorem}}
\newcommand{\ec}{\end{corollary}}
\newcommand{\el}{\end{lemma}}
\newcommand{\ee}{\end{example}}
\newcommand{\ep}{\end{proposition}}
\newcommand{\eq}{\end{question}}
\newcommand{\ed}{\end{definition}}
\newcommand{\epr}{{\ \vbox{\hrule\hbox{%
\vrule height1.3ex\hskip0.8ex\vrule}\hrule}}\\\par}

\def\R{\mathbb{R}}
\def\Z{\mathbb{Z}}
\def\supp{\mbox{supp}}
\def\rank{\mbox{rank}}
\def\x{\hat{x}}
\def\td{{\rm td}}

\begin{document}

\title{\bf Circuit and Graver Walks\\and\\Linear and Integer Programming}

\author{
Shmuel Onn\thanks{\small Technion - Israel Institute of Technology. Email: onn@technion.ac.il}
}

\date{}

\maketitle

\begin{abstract}
We show that a circuit walk from a given feasible point of a given linear program
to an optimal point can be computed in polynomial time using only linear algebra
operations and the solution of the single given linear program.

We also show that a Graver walk from a given feasible point of a given integer program
to an optimal point is polynomial time computable using an integer programming
oracle, but without such an oracle, it is hard to compute such a walk even if
an optimal solution to the given program is given as well.

Combining our oracle algorithm with recent results on sparse integer programming,
we also show that Graver walks from any point are polynomial time computable
over matrices of bounded tree-depth and subdeterminants.

\end{abstract}

\section{Introduction}

Walks from a given feasible point of a given linear program or integer program to an optimal point
had been studied for decades, ever since the introduction of the Simplex method and
the Hirsch conjecture. Here we show the following. For linear programming, we show
that a {\em circuit walk} can be computed in polynomial time by solving a single linear
program. For integer programming, we show
that it is hard to compute a {\em Graver walk} by solving a single integer program,
but such a walk can be computed in polynomial time using an integer programming oracle.
We conclude with several remaining open problems. Following is a more detailed description.

\subsection{Circuit walks and linear programming}

We consider the following linear programming problem and corresponding polytope,
$$\mbox{LP}:\ \max\{wx\, :\, x\in\R^n,\ Ax=b,\, l\leq x\leq u\},
\ \ \mbox{P}\, :=\, \{x\in\R^n\, :\, Ax=b,\, l\leq x\leq u\},$$
where $A$ is an integer $m\times n$ matrix and $b\in\Z^m$ and $w,l,u\in\Z^n$ are integer vectors.
For simplicity we assume that all data is integer, but this results in no loss of generality
since, given rational data, we can always scale it, and make a change of variables by
scaling, to obtain a polynomially equivalent problem with integer data.

A {\em walk} is a sequence $x^0,x^1,\dots,x^s=x^*$ of points feasible in LP, of increasing
objective value, starting at a given point $x^0$ and terminating at point $x^*$ optimal for LP.
(We consider only such {\em monotone} walks and for brevity omit the term monotone.)

An {\em edge walk} is one where each $x^k$ is a vertex of $P$, which is obtained from $x^{k-1}$ by
taking a maximal step along an edge of $P$. A major open problem in the theory of convex
polytopes is whether, starting at any given vertex $x^0$, an edge walk with a polynomial
number $s$ of steps always exists. A stronger, major open problem in linear
programming, is whether an edge walk could be computed in polynomial time.

It is natural to relax these hard long open problems by considering a more flexible type of
walks. A {\em circuit walk} is one with each $x^k$ obtained from $x^{k-1}$ by
taking a maximal step along a {\em circuit} of $A$ (see Section 2 for the precise definitions).
Each edge walk is in particular a circuit walk. Circuit walks were first introduced and studied
in \cite{BDF,BFH}, in the hope that they might eventually lead to a strongly polynomial time
algorithm for linear programming, a fundamental open problem in the theory of linear programming;
see these references for more details and motivation for studying such walks.
The two analog problems for circuit walks are
whether a circuit walk from any point with a polynomial number of steps always exists, and whether
a circuit walk could be computed in polynomial time. For the first problem, it has been known
for quite some time that the answer is positive and a circuit walk with a polynomial number
of steps {\em does} always exist; in fact, such a short circuit walk is obtained when each step
is taken along a circuit giving the maximum possible objective improvement \cite{DHL,HOW}.
So it is tempting to ask whether such best circuit steps could be computed in polynomial time.
Unfortunately, it was shown in \cite{DKS} that finding a best circuit step is NP-hard, leaving
open the problem of whether it is possible to compute a circuit walk in polynomial time by other
means. Here we answer this affirmatively. Moreover, our algorithm uses only linear algebra
operations along with the solution of a single linear program, the given one.

\bt{main}
A circuit walk from any point is polynomial time computable using only linear
algebra operations along with the solution of a single linear program.
\et

While, as noted, finding a best circuit step is NP-hard, in proving this theorem we
show that finding a {\em good enough} circuit step {\em can} be done in polynomial time.

\subsection{Graver walks and integer programming}

We also consider the following integer programming counterpart of our problem,
$$\mbox{IP}:\ \max\{wx\, :\, x\in\Z^n,\ Ax=b,\, l\leq x\leq u\}\ .$$
A {\em Graver walk} is one where each $x^k$ is feasible in IP and obtained from $x^{k-1}$ by taking
a maximal step along a {\em Graver-element} of $A$ (see Section 3 for the precise definitions).\break
The two analog problems for Graver walks are whether a Graver walk from any integer
point with a polynomial number of steps always exists, and whether a Graver walk could be
computed in polynomial time. For the first problem, it has been known for quite some time
that the answer is positive and a Graver walk with a polynomial number of steps {\em does}
always exist; in fact, such a short Graver walk is obtained when each step is taken along a
Graver-element giving the maximum possible objective improvement \cite{HOW,Onn}. This implies that
integer programming, which is NP-hard, reduces to finding a best Graver step, and so unfortunately,
finding a best Graver step is NP-hard as well. But it is tempting to ask whether an integer analog of Theorem \ref{main}
holds, namely, given an initial point $x^0$ feasible for IP {\em and an optimal solution} $\x$
for IP, can we compute a Graver walk in polynomial time? Unfortunately the integer analogs
of Theorem \ref{main} and its proof fail, and on the negative side we show, in \S 3.1,
that computing a Graver walk is NP-hard even given an optimal solution $\x$ for IP.
But on the positive side, in \S 3.2, we show that, given access to an integer programming oracle,
that is, by solving multiple integer programs, we {\em can} compute a Graver walk in polynomial
time. As usual, $\Delta(A)$ is the maximum absolute value of any subdeterminant of $A$.

\bt{integer}
A Graver walk from any point is computable in time which is polynomial in the
bit size of $w,l,u$ and in $\Delta(A)$ using an integer programming oracle.
\et

Combining this theorem with recent results on sparse integer programming from \cite{EHKKLO,KLO,KO},
we also show in Theorem \ref{tree-depth} that Graver walks can
be computed in polynomial time over matrices of bounded tree-depth and subdeterminants.

\subsection{Notation}

The {\em support} of $x\in\R^n$ is the set $\supp(x):=\{i:x_i\neq 0\}\subseteq[n]$ of indices of
its nonzero entries. A {\em circuit} of $A$ is a nonzero kernel element with minimal support,
that is, a vector $c\in\Z^n$, $c\neq 0$, $Ac=0$, such that there is no $d\in\R^n$, $d\neq 0$,
$Ad=0$, with $\supp(d)\subsetneqq\supp(c)$. A vector $c$ is {\em dominated} by a vector $h$
if $\supp(c)\subseteq\supp(h)$ and $c_ih_i\geq 0$ for all $i$. More strongly, $c$ is
{\em conformal} to $h$, denoted $c\preceq h$, if $c_ih_i\geq 0$ and $|c_i|\leq|h_i|$ for all $i$.
A {\em Graver-element} of $A$ is a nonzero integer kernel element which is conformal-minimal,
that is, a vector $g\in\Z^n$, $g\neq 0$, $Ag=0$, such that there is no $f\in\Z^n$, $f\neq 0$,
$Af=0$, with $f\precneqq g$. The $j$th column of $A$ is $A^j$ and for $J\subseteq[n]$ we let
$A^J:=\{A^j:j\in J\}$. The {\em matroid} of $A$ is the matroid on $[n]$ of linear dependencies
on $A$, where a subset $C\subseteq[n]$ is a matroid circuit if and only if $A^C$ is a linearly
dependent set and $A^{C\setminus\{j\}}$ is a linearly independent set for each $j\in C$.

\section{Circuit walks and linear programming}

A {\em circuit walk} is a sequence $x^0,x^1,\dots,x^s=x^*$ of points feasible for LP,
that is, points in $P$, terminating at a point $x^*$ optimal for LP, such that for
all $k$ we have that $wx^k>wx^{k-1}$ and $c^k:=x^k-x^{k-1}$ is a circuit $c^k$ of $A$ forming
a maximal step in direction $c^k$ at $x^{k-1}$, that is, $x^{k-1}+\alpha c^k\notin P$ for
all real $\alpha>1$. All non initial points in a circuit walk are on the boundary of $P$
but are generally not vertices.

The next lemma is known in one form or another, see e.g. \cite[Lemma 2.17]{Onn}
or \cite{BV}. For completeness and convenience of the reader we include a short constructive
proof followed by a detailed example illustrating the intricate underlying algorithm.

\bl{sum}
Given any $h\neq0$ in the kernel of $A$ we can find in polynomial time with only linear
algebra operations $t\leq n$ circuits $c^k$ of $A$ conformal to $h$ with $h=\sum_{k=1}^tc^k$.
\el
\bpr
We start with a sequence of procedures used repeatedly in the algorithm.

\vskip.2cm
{\bf Finding a matroid circuit:} Given $H\subseteq[n]$ with $\rank(A^H)<|H|$
find a matroid circuit $C\subseteq H$ as follows. Initialize $J':=\emptyset$ and
$J:=\{j\}$ for some $j\in H$. While $\rank(A^J)>\rank(A^{J'})$ do:
set $J':=J$; set $J:=J'\uplus\{j\}$ for some $j\in H\setminus J'$; repeat.
Return the matroid circuit $C:=\{j\in J:\rank(A^{J\setminus\{j\}})=\rank(A^J)\}$.

\vskip.2cm
{\bf Reducing by a circuit:} Given $h\neq0$ in the kernel of $A$ and a circuit $c$
of $A$ with $\supp(c)\subseteq\supp(h)$, {\em the reduction of $h$ by $c$} is
$h':=h-d$ for a suitable and easily computable scalar multiple $d$ of $c$ such that
$h'$ is dominated by $h$ and has at least one more zero entry than $h$,
that is, $h'_ih_i\geq 0$ for all $i$ and $\supp(h')\subsetneqq\supp(h)$.

\vskip.2cm
{\bf Finding a dominated circuit:} Let $h\neq0$ be a given vector in the kernel of $A$.
Set $k:=0$ and $g^0:=h$. While $g^k\neq 0$ do: find a matroid circuit
$C^{k+1}\subseteq\supp(g^k)$; let $c^{k+1}$ be the unique up to scalar multiple circuit
of $A$ with $\supp(c^{k+1})=C^{k+1}$; let $g^{k+1}$ be the reduction of $g^k$ by $c^{k+1}$;
set $k:=k+1$; repeat. Let $h=g^0,\dots,g^r=0$ be the sequence obtained. The process is
repeated $r\leq n$ times since $\supp(g^{k+1})\subsetneqq\supp(g^k)$ for all $k$. All $g^k$
are dominated by $h$, and the dominated $g^{r-1}$ is a multiple of $c^r$ hence is a
circuit of $A$. Note that if any circuit $c^k$ obtained on the way or its negation $-c^k$
happens to be dominated by $h$ then we can use it and stop the process earlier.

\vskip.2cm
We can now describe the algorithm proving the lemma. Let $h\neq0$ be a given vector in the
kernel of $A$. Set $k:=0$ and $h^0:=h$. While $h^k\neq 0$ do: find a circuit $c^{k+1}$ dominated
by $h^k$; let $h^{k+1}:=h^k-d^{k+1}$ be the reduction of $h^k$ by $c^{k+1}$ and note that $d^{k+1}$
and $h^{k+1}$ must be conformal to $h^k$ and hence also to $h$; set $k:=k+1$;
repeat. Let $h=h^0,\dots,h^t=0$ be the sequence obtained. The process is repeated at most
$t\leq n$ times since $\supp(h^{k+1})\subsetneqq\supp(h^k)$ for all $k$, and $h=\sum_{k=1}^t d^k$
is the sum of the $t\leq n$ obtained circuits $d^k$ of $A$ which are conformal to $h$ as desired.
\epr

\be{example1}
{\rm
We demonstrate the algorithm of the lemma on the following input,
$$A\ =\
\begin{pmatrix}
2&1&0&1&2&2&1\\
0&1&1&1&0&1&0\\
2&2&1&0&0&1&1\\
\end{pmatrix},\quad
h^0\ :=\ h\ =\
\begin{pmatrix}
2&5&-9&1&-4&3&-8
\end{pmatrix}\ .$$

We let $g^0:=h^0$, find a circuit
$c^1:=\begin{pmatrix}-1&2&-2&0&0&0&0\end{pmatrix}$, reduce $g^0$ by it to
$g^1:=\begin{pmatrix}0&9&-13&1&-4&3&-8\end{pmatrix}$,
find a circuit
$c^2:=\begin{pmatrix}0&-1&2&-1&1&0&0\end{pmatrix}$, observe that $-c^2$ is already
dominated by $h^0$, stop the process and use this circuit.

Now we reduce $h^0$ by the dominated circuit just found and obtain
$$h^1\ :=\ h^0-d^1\ =\ \begin{pmatrix}2&4&-7&0&-3&3&-8\end{pmatrix},
\quad d^1\ :=\ \begin{pmatrix}0&1&-2&1&-1&0&0\end{pmatrix}\ .$$

We now let $g^0:=h^1$, find a circuit $c^1:=\begin{pmatrix}-1&2&-2&0&0&0&0\end{pmatrix}$,
reduce $g^0$ by it to $g^1:=\begin{pmatrix}0&8&-11&0&-3&3&-8\end{pmatrix}$,
find a circuit $c^2:=\begin{pmatrix}0&0&-1&0&-1&1&0\end{pmatrix}$, observe that $c^2$ is already
dominated by $h^1$, stop the process and use this circuit.

Now we reduce $h^1$ by the dominated circuit just found and obtain
$$h^2\ :=\ h^1-d^2\ =\ \begin{pmatrix}2&4&-4&0&0&0&-8\end{pmatrix},
\quad d^2\ :=\ \begin{pmatrix}0&0&-3&0&-3&3&0\end{pmatrix}\ .$$

We now let $g^0:=h^2$, find a circuit $c^1:=\begin{pmatrix}-1&2&-2&0&0&0&0\end{pmatrix}$,
reduce $g^0$ by it to $g^1:=\begin{pmatrix}0&8&-8&0&0&0&-8\end{pmatrix}$,
find a circuit $c^2:=\begin{pmatrix}0&1&-1&0&0&0&-1\end{pmatrix}$, reduce $g^1$ by it to $g^2:=0$,
and the process ends with $g^1$ (or $c^2$) dominated by $h^2$.

Now we reduce $h^2$ by the dominated circuit just found and obtain
$$h^3\ :=\ h^2-d^3\ =\ \begin{pmatrix}2&0&0&0&0&0&-4\end{pmatrix},
\quad d^3\ :=\ \begin{pmatrix}0&4&-4&0&0&0&-4\end{pmatrix}\ .$$

We now let $g^0:=h^3$, find a circuit $c^1:=\begin{pmatrix}1&0&0&0&0&0&-2\end{pmatrix}$,
observe that it is a multiple of $h^3$ and hence dominated by $h^3$,
so we can reduce $h^3$ by $c^1$ to obtain
$$h^4\ :=\ h^3-d^4\ =\ 0,\quad d^4\ :=\ \begin{pmatrix}2&0&0&0&0&0&-4\end{pmatrix}\ .$$

So we found the desired expression $h=\sum_{k=1}^4d^k$ as a sum of conformal circuits.
}
\ee

\vskip.2cm\noindent
{\bf Theorem \ref{main}} A circuit walk from any point is polynomial time computable
using only linear algebra operations along with the solution of a single linear program.

\vskip.2cm
\bpr
Let $x^0$ be the given feasible point of LP. If $A$ is a zero matrix then $b={\bf 0}$
and an optimal solution to LP is $x^*$ with $x^*_i=u_i$ if $w_i\geq0$ and
$x^*_i=l_i$ if $w_i<0$, and a short circuit walk is trivially obtained by replacing
one by one each entry of $x^0$ by the corresponding entry of $x^*$. So we may assume
$\Delta:=\Delta(A)\geq 1$. Use linear programming to find in polynomial time
an optimal solution $\x$ of LP. Let $h^0:=\x-x^0$, which is in the kernel of $A$.
Assume $h^0\neq 0$ else $x^0=\x$ is already optimal. Use Lemma \ref{sum} to find in
polynomial time $t\leq n$ circuits $c^i$ of $A$
conformal to $h^0$ such that $h^0=\sum_{i=1}^t c^i$. Since $l\leq x^0\leq u$
and $l\leq \x=x^0+\sum_{\i=1}^t c^i\leq u$ and each $c^i$ is conformal to $h^0$,
we have that for all $i$ also $l\leq x^0+c^i\leq u$ and hence $x^0+c^i$
is in $P$. As $w\x-wx^0=\sum_{i=1}^t wc^i$, some $c^k$ satisfies
$$wc^k\ \geq\ {1\over t}(w\x-wx^0)\ \geq\ {1\over n}(w\x-wx^0)\ .$$
Let $d^k$ be the maximal positive multiple of $c^k$ such that $x^1:=x^0+d^k$ is in $P$,
easily computed in polynomial time by checking the lower and upper bounds $l$ and $u$. Then
$$wx^1-wx^0\ =\ wd^k\ \geq\ wc^k\ \geq\ {1\over n}(w\x-wx^0)$$
which implies
$$w\x-wx^1\ \leq\ \left(1-{1\over n}\right)(w\x-wx^0)\ .$$
Now let $h^1:=\x-x^1$ and repeat the process, obtaining a new point $x^2$ in $P$ with
$$w\x-wx^2\ \leq\ \left(1-{1\over n}\right)(w\x-wx^1)
\ \leq\ \left(1-{1\over n}\right)^2(w\x-wx^0)\ .$$
Repeating this we obtain a sequence of points $x^0,x^1,x^2,\dots$ in $P$ where
$$w\x-wx^k\ \leq\ \left(1-{1\over n}\right)^k(w\x-wx^0)\ .$$
Let $r:=n\lceil\ln(\Delta(w\x-wx^0))\rceil$. Taking the natural logarithm on both sides
of the above expression for $k:=r$ and using the Taylor bound $\ln(1-{1\over n})<-{1\over n}$
we obtain
\begin{eqnarray*}
\ln(w\x-wx^r) & \leq & r\cdot\ln\left(1-{1\over n}\right)+\ln(w\x-wx^0)
 \ < \ -r\cdot{1\over n}+\ln(w\x-wx^0)             \\
& \leq & -n\cdot\ln(\Delta(w\x-wx^0)){1\over n}+\ln(w\x-wx^0)\ =\ \ln\left({1\over\Delta}\right)\ ,
\end{eqnarray*}
which implies $w\x-wx^r<{1\over\Delta}$.
Since $\Delta\leq\max|A_{i,j}|^mm^{{m\over2}}$, the number $r$ of points needed
is polynomial in the bit size of $A,w,l,u$. 
Now, as is well known, any vertex is a basic solution, hence by Cramer's rule,
its components are integer multiples of ${1\over|\Gamma|}\geq{1\over\Delta}$ where
$\Gamma$ is the determinant of some square submatrix of $A$, and since $w$ is integer,
so is its objective value. Therefore the gap between the objective values of optimal vertices
and non optimal vertices is bounded below by $1\over\Delta$, and therefore
each face of $P$ containing $x^r$ must contain some optimal vertex, in particular the
smallest such face $F$ which, with $L:=\{i:x^r_i=l_i\}$ and $U:=\{i:x^r_i=u_i\}$, is
$$F\ =\ \{x\in\R^n\ :\ Ax=b,\ x_i=l_i,\ i\in L,
\ x_i=u_i,\ i\in U,\ l_i\leq x_i\leq u_i,\ i\notin L\uplus U\}\ .$$
Let $H:=[n]\setminus(L\uplus U)$. If $\rank(A^H)=|H|$ then $x^r$ is a vertex and
we are done. Otherwise we find a matroid circuit $C\subseteq H$, and let $c$ be the unique
circuit with $\supp(c)=C$ such that $wc\geq0$ and $c$ forms a maximal step in direction $c$,
namely, the point $x^{r+1}:=x^r+c$ is in $P$ hence satisfies $l\leq x^{r+1}\leq u$,
and $x^{r+1}_i=l_i$ or $x^{r+1}_i=u_i$ for at least one $i\notin L\uplus U$.
So $x^{r+1}$ is on the boundary of $F$ and hence on a lower dimensional face of
$P$ which again must contain some optimal vertex since $wx^{r+1}\geq wx^r$.
Repeating this, we take less than $n$ maximal circuit steps, each one in turn to a lower
dimensional face, till we reach an optimal vertex $x^s$ with $s<r+n$. The sequence
$x^0,\dots,x^r,\dots,x^s=:x^*$ is a polynomially computable circuit walk.
\epr

\be{example2}
{\rm
We demonstrate the algorithm of the theorem on the following input,
$$A\ =\
\begin{pmatrix}
2&1&0&1&2&2&1\\
0&1&1&1&0&1&0\\
2&2&1&0&0&1&1\\
\end{pmatrix},\quad
b\ =\
\begin{pmatrix}
0\\
0\\
0\\
\end{pmatrix}\ ,$$
$$
l\ =\
\begin{pmatrix}
-5&-5&0&-5&0&-5&0
\end{pmatrix},\quad
u\ =\
\begin{pmatrix}
0&0&9&0&9&0&9
\end{pmatrix}\ ,$$
$$w\ =\
\begin{pmatrix}
1&1&-1&1&-1&1&-1,
\end{pmatrix},\quad
x^0\ =\
\begin{pmatrix}
-2&-5&9&-1&4&-3&8
\end{pmatrix}\ .$$
Linear programming gives an optimal solution $\x=0$, and $\Delta=4$ gives a lower bound
$1\over4$ on the gap between the objective values of optimal and non optimal vertices.

We now let $h^0:=\x-x^0$ which is a conformal sum $h^0=\sum_{k=1}^4 d^k$
of circuits as in Example \ref{example1}. We find the maximum $wd^k$ which
is $wd^3=12>{32\over 7}={1\over n}(w\x-wx^0)$.

We take a maximal step along $d^3$ taking us to $x^1:=x^0+{5\over4}d^3$.
We continue with $h^1:=\x-x^1=\begin{pmatrix}2&0&-4&1&-4&3&-3\end{pmatrix}$.
We use the lemma to write it as a conformal sum of circuits $h^1=\sum_{k=1}^3c^k$ with
$c^1=\begin{pmatrix}{1\over 2}&0&-1&1&-1&0&0\end{pmatrix}$,
$c^2=\begin{pmatrix}0&0&-3&0&-3&3&0\end{pmatrix}$,
and $c^3=\begin{pmatrix}{3\over 2}&0&0&0&0&0&-3\end{pmatrix}$.
We compute the values $wc^k$ and the maximum is $wc^2=9$.
We take a maximal step along $c^2$ taking us to $x^2:=x^1+c^2$.
We continue with $h^2:=\x-x^2=\begin{pmatrix}2&0&-1&1&-1&0&-3\end{pmatrix}$.
We see (or use the lemma) that it is the conformal sum of circuits $h^2=c^1+c^3$.
We compute the values $wc^k$ and the maximum is $wc^3={9\over2}$.
We take a maximal step along $c^3$ taking us to $x^3:=x^2+c^3$.
We continue with $h^3:=\x-x^3=\begin{pmatrix}{1\over2}&0&-1&1&-1&0&0\end{pmatrix}$.
We see that $h^3=c^1$ is itself a circuit, so we take a maximal step along $c^1$ taking
us to $x^4:=x^3+c^1=0=\x$ which is already an optimal vertex. So here we arrive
directly to an optimal vertex without the second part of the proof of going down to lower
and lower dimensional faces. The resulting circuit walk is $x^0,x^1,x^2,x^3,x^4=x^*$ with
$$x^0\ =\
\begin{pmatrix}
-2&-5&9&-1&4&-3&8
\end{pmatrix},\quad
x^1\ =\
\begin{pmatrix}
-2&0&4&-1&4&-3&3
\end{pmatrix}
\ ,$$
$$x^2\ =\
\begin{pmatrix}
-2&0&1&-1&1&0&3
\end{pmatrix},\ \
x^3\ =\
\begin{pmatrix}
-{1\over2}&0&1&-1&1&0&0
\end{pmatrix},\ \
x^4\ =\ x^*\ =\ 0\ .
$$
}
\ee

\section{Graver walks and integer programming}

A {\em Graver walk} is a sequence $x^0,x^1,\dots,x^s=x^*$ of integer points in $P$
terminating at $x^*$ optimal for IP, such that for all $k$ we have that $wx^k>wx^{k-1}$
and $f^k:=x^k-x^{k-1}$ is a positive integer multiple of some Graver-element $g^k$
of $A$ forming a maximal step in direction $g^k$ at $x^{k-1}$, that is,
$x^{k-1}+\alpha f^k$ is not in $P\cap\Z^n$ for all real $\alpha>1$.

\subsection{Impossibility}

We begin by showing that deciding non Graver-elements is NP-complete.

\bp{Graver}
It is NP-complete to decide if a given integer vector $h\in\Z^n$ is not a Graver-element
of a given integer ${m\times n}$ matrix $A$, even when $A$ has a single row.
\ep
\bpr
The {\em subset-sum problem}, well known to be NP-complete, is to decide, given positive
integers $a_0,a_1,\dots,a_r$, whether there is a subset $I\subseteq[r]$ with $\sum_{i\in I}a_i=a_0$.

Given such $a_i$, we let $m:=1$, $n:=r+2$, $A:=[a_1,\dots,a_r,-a_0,a_0-\sum_{i=1}^ra_i]$, and
$h:={\bf 1}:=[1,\dots,1]$ the all-ones vector in $\Z^n$, which is clearly in the kernel of $A$.

Suppose there is an $I\subseteq[r]$ with $\sum_{i\in I}a_i=a_0$. Define $f\in\{0,1\}^n$
by $f_i:=1$ if $i\in I$ or $i=r+1$, and $f_i:=0$ otherwise. Then $f$ is integer, $f\neq0$,
$Af=\sum_{i\in I}a_i-a_0=0$, and $f\precneqq h$, hence $h$ is not a Graver-element.
Conversely, suppose $h$ is not a Graver-element. Then there exists an integer $f\neq0$,
$Af=0$, and $f\precneqq h$. Since $f$ is conformal to $h$, we have $f\in\{0,1\}^n$.
Consider $f_{r+1}$ and $f_{r+2}$. If $f_{r+1}=f_{r+2}=0$ then $f\neq0$ implies $f_i=1$
for some $i\in[r]$ but then $0=Af\geq a_i>0$ which is impossible. If $f_{r+1}=f_{r+2}=1$ then
$f\neq h$ implies $f_i=0$ for some $i\in[r]$ but then $0=Af\leq -a_i<0$ which is impossible.
If $f_{r+1}=1$ and $f_{r+2}=0$ then letting $I:=\{i\in[r]:f_i=1\}$ we have
$0=Af=\sum_{i\in I}a_i-a_0$ hence $\sum_{i\in I}a_i=a_0$.
Finally, if $f_{r+1}=0$ and $f_{r+2}=1$ then letting $I:=\{i\in[r]:f_i=0\}$ we have
$0=Af=\sum_{i\in [r]\setminus I}a_i+a_0-\sum_{i=1}^ra_i$ hence again $\sum_{i\in I}a_i=a_0$.
This shows that $h$ is not a Graver-element if and only if
there is some subset $I\subseteq[r]$ with $\sum_{i\in I}a_i=a_0$.
\epr

The integer Carath\'eodory theorem asserts that for every nonzero integer kernel element
$h$ of $A$ there exist $t\leq 2n-2$ Graver-elements $g^k$ of $A$ conformal to $h$
and positive integers $\lambda_k$ such that $h=\sum_{k=1}^t\lambda_k g^k$, see \cite{Onn}.
It is tempting to ask whether an integer analog of Lemma \ref{sum} holds, namely,
given $h$, can we compute such $\lambda_k$ and $g^k$ in polynomial time?
The next proposition shows that the answer is negative.

\bp{conformal}
It is NP-hard, given an integer kernel element $h\neq 0$ of given $A$, to find Graver-elements
$g^k$ conformal to $h$ and positive integers $\lambda_k$ with $h=\sum_{k=1}^t\lambda_k g^k$.
\ep
\bpr
If $h$ is a Graver-element then the only possible such expression is the trivial sum $h=h$,
that is, $h=\sum_{k=1}^t\lambda_k g^k$ with $t=1$, $\lambda_1=1$, $g^1=h$. If $h$ is not
a Graver-element then any such expression must either involve some coefficient $\lambda_k\geq 2$
or involve $t\geq 2$ summands. So if we could find such an expression then we could decide if
$h$ is a Graver-element which is NP-complete by Proposition \ref{Graver} above.
\epr

Finally, it is tempting to ask whether an integer analogs of Theorem \ref{main} holds,
namely, given an initial point $x^0$ feasible for IP {\em and an optimal solution} $\x$
for IP, can we compute a Graver walk in polynomial time? The answer is again negative.

\bp{Graver_walk}
Given a feasible point $x^0$ to a given integer programming problem
$$\mbox{IP}:\ \max\{wx\, :\, x\in\Z^n,\ Ax=b,\, l\leq x\leq u\}\ ,$$
it is NP-hard to compute a Graver walk, even when $A$ has a single row, $b=0$, $l={\bf 0}$
is all-zeros, $w=u={\bf 1}$ are all-ones, and an optimal solution $\x$ for IP is given too.
\ep
\bpr
As in the proof of Proposition \ref{Graver}, given input to the subset-sum problem, let $m:=1$,
$n:=r+2$, and $A:=[a_1,\dots,a_r,-a_0,a_0-\sum_{i=1}^ra_i]$. Also let $b:=0$, $x^0:=l:={\bf 0}$,
$\x:=w:=u:={\bf 1}$. So $x^0$ is feasible and $\x$ is the only optimal solution.

Suppose $x^0,\dots,x^s=x^*=\x$ is a Graver walk. If it has one step namely $s=1$ then
${\bf 1}=x^1-x^0$ is a Graver-element, hence as shown in the proof of Proposition \ref{Graver},
there is no $I\subseteq[r]$ with $\sum_{i\in I}a_i=a_0$. If it has more than one
step namely $s\geq 2$ then $0\neq x^1-x^0=x^1\precneqq{\bf 1}$ hence ${\bf 1}$ is not a
Graver-element, and hence as shown in the proof of Proposition \ref{Graver}, there is some
$I\subseteq[r]$ with $\sum_{i\in I}a_i=a_0$. Therefore, if we could compute a Graver walk
then we could decide the subset-sum problem.
\epr

\subsection{An oracle algorithm}

In contrast to the impossibility result of Proposition \ref{Graver_walk}, we now show that,
given access to an integer programming oracle, that is, by solving multiple integer programs,
we {\em can} compute a Graver walk in polynomial time. We begin with a lemma.
\bl{sub_Graver}
Given any $h\neq0$ in the kernel of any given $A$, a Graver-element of $A$ which
is conformal to $h$ can be obtained by solving an integer programming problem.
\el
\bpr
Let $\sigma_i:=1$ if $h_i\geq0$ and $\sigma_i:=-1$ if $h_i<0$. Consider the following program,
\begin{equation}\label{auxiliary}
\!\!\min\left\{\|x\|_1=\sum_{i=1}^n\sigma_ix_i\, :\,  x\in\Z^n,\ Ax=0,\,
\sum_{i=1}^n\sigma_ix_i\geq 1,\, 0\leq \sigma_ix_i\leq \sigma_ih_i,\, i=1,\dots,n\right\}\, .
\end{equation}
Any optimal solution $x$ satisfies $x\preceq h$ by its feasibility and there is no nonzero kernel
element $g$ with $g\precneqq x$ by its optimality, hence $x$ is a Graver-element conformal to $h$.
\epr

\vskip.2cm\noindent{\bf Theorem \ref{integer}}
A Graver walk from any point is computable in time which is polynomial in the
bit size of $w,l,u$ and in $\Delta(A)$ using an integer programming oracle.

\vskip.2cm
\bpr
Let $x^0$ be the given feasible point of IP. If $A$ is a zero matrix then $b={\bf 0}$
and an optimal solution to IP is $x^*$ with $x^*_i=u_i$ if $w_i\geq0$ and
$x^*_i=l_i$ if $w_i<0$, and a short Graver walk is trivially obtained by replacing
one by one each entry of $x^0$ by the corresponding entry of $x^*$. So we may assume
$\Delta:=\Delta(A)\geq 1$. Use the oracle to find an optimal point $\x$ of IP.
Let $h^0:=\x-x^0$ and assume $h^0\neq 0$ since otherwise $x^0=\x$ is already optimal.
This in particular implies $\rank(A)+1\leq n$. Use Lemma \ref{sum} to write
$h^0=\sum_{i=1}^t c^i$ as a conformal sum of $t\leq n$ circuits $c^i$ of $A$. For each $i$ let
$\alpha_i$ be the positive rational number such that the entries of $d^i:={1\over\alpha_i}c^i$
are relatively prime integers, that is, $d^i\in\Z^n$ and $\gcd(d^i_1,\dots,d^i_n)=1$.
Then each $d^i$ satisfies $\|d\|_1\leq(\rank(A)+1)\Delta(A)$, see \cite[Lemma 3.18]{Onn},
and $d^i$ is a Graver-element. Let $f^0:=\sum_{i=1}^t(\alpha_i-\lfloor\alpha_i\rfloor)d^i$.
As the sum is conformal and $\alpha_i-\lfloor\alpha_i\rfloor<1$ for all $i$,
$$\|f^0\|_1\ =\ \sum_{i=1}^t(\alpha_i-\lfloor\alpha_i\rfloor)\|d^i\|_1
\ <\ \sum_{i=1}^t\|d^i\|_1\ \leq\ tn\Delta\ \leq\ n^2\Delta\ .$$
Since $f^0=h^0-\sum_{i=1}^t\lfloor\alpha_i\rfloor d^i$ we have that $f^0$ is an integer vector
in the kernel of $A$. So we can apply Lemma \ref{sub_Graver} and obtain a Graver-element $g^1$
such that $g^1\preceq f^0$. We can now apply Lemma \ref{sub_Graver} again and obtain a
Graver-element $g^2$ such that $g^2\preceq f^0-g^1$.
Continuing to apply Lemma \ref{sub_Graver} repeatedly this way, we can obtain a sequence
$g^1,g^2,\dots$ of Graver-elements such that $g^k\preceq f^0-\sum_{i<k}g^i$ for all $k$.
Since $\|f^0\|_1\leq n^2\Delta$ and $\|g^i\|_1\geq 1$ for all $i$, the sequence ends at
some $g^r$ with $r\leq n^2\Delta$ and $f^0=\sum_{i=1}^r g^i$, and we get an expression
$$h^0\ =\ f^0+\sum_{i=1}^t\lfloor\alpha_i\rfloor d^i\ =\
\sum_{i=1}^{r+t}\lambda_i g^i,\quad \lambda_i:=1,\ i\leq r,\quad
g^i:=d^{i-r},\ \lambda_i:=\lfloor\alpha_{i-r}\rfloor,\ i>r$$
of $h^0$ as a positive integer combination of $r+t\leq 2n^2\Delta$ conformal Graver-elements.
As in the proof of Theorem \ref{main}, there is some $g^k$ with
$w\lambda_kg^k\ \geq\ {1\over 2n^2\Delta}(w\x-wx^0)$.
Determine the largest positive integer $\lambda\geq\lambda_k$ such that
$x^1:=x^0+\lambda g^k$ satisfies the lower and upper bounds, which implies
$w\x-wx^1\leq(1-{1\over{2n^2\Delta}})(w\x-wx^0)$. Now let $h^1:=\x-x^1$
and repeat the process, obtaining a new point $x^2$, and so on, obtaining
a sequence of points $x^0,x^1,x^2,\dots$ which are feasible in IP, where for all $k$ we have
$$w\x-wx^k\ \leq\ \left(1-{1\over{2n^2\Delta}}\right)^k(w\x-wx^0)\ .$$
Taking the natural logarithm on both sides and using the Taylor bound as in the proof
of Theorem \ref{main}, we find that, with $s:=2n^2\Delta\lceil\ln(w\x-wx^0)\rceil$,
we have $w\x-wx^s<1$, and since $w,\x,x^s$ are integer vectors we must in fact have $wx^s=w\x$
and hence $x^s$ is optimal too. Thus, the resulting sequence $x^0,\dots,x^s=:x^*$ is a
Graver walk computed in time polynomial in the bit size of $w,l,u$ and in $\Delta=\Delta(A)$.
\epr

\subsection{Bounded tree-depth and subdeterminants}

Recent results on sparse integer programming, see \cite{EHKKLO,KLO,KO}, imply that for
any fixed $a,d$, the integer programming problem IP can be solved in polynomial
time over any matrix $A$ satisfying $\max|A_{i,j}|\leq a$ and $\min\{\td(A),\td(A^T)\}\leq d$,
where $\td(A)$ is the {\em tree-depth} of $A$, a parameter which plays a central role in sparsity,
see \cite{NO}, and which is defined as follows. The {\em height} of a rooted tree is the maximum
number of vertices on a path from the root to a leaf. Given a graph $G=(V,E)$, a rooted tree on
$V$ is {\em valid} for $G$, if for each edge $\{j,k\}\in E$, either $j$ lies on the path from the
root to $k$ or vice versa. The {\em tree-depth} $\td(G)$ of $G$ is the smallest height of a rooted
tree which is valid for $G$. The graph of an $m\times n$ matrix $A$ is the graph $G(A)$ on $[n]$
where ${j,k}$ is an edge if and only if there is an $i\in[m]$ such that $A_{i,j}A_{i,k}\neq 0$.
The {\em tree-depth} of the matrix $A$ is the tree-depth $\td(A):=\td(G(A))$ of its graph.

Recent results on sparse integer programming from \cite{EHKKLO,KLO,KO} imply the following.

\bp{sparse_IP} Fix any $a,d$. The integer program IP over any matrix $A$ satisfying
$\max|A_{i,j}|\leq a$ and $\min\{\td(A),\td(A^T)\}\leq d$ can be solved in polynomial time.
\ep

\bt{tree-depth}
Fix any $a,d$. A Graver walk from any point is polynomial time computable over
any matrix $A$ which satisfies $\Delta(A)\leq a$ and $\min\{\td(A),\td(A^T)\}\leq d$.
\et
\bpr
We apply the oracle algorithm of Theorem \ref{integer}. The algorithm solves IP once to
obtain $\x$, and then repeatedly solves an auxiliary integer program \eqref{auxiliary} as
in the proof of Lemma \ref{sub_Graver}. IP is in standard form with matrix $A$,
but \eqref{auxiliary} has the additional inequality $\sum_{i=1}^n\sigma_ix_i\geq 1$.
Instead of solving \eqref{auxiliary}, we can solve, for $k=1,\dots,n$,
$$\min\left\{\|x\|_1=\sum_{i=1}^n\sigma_ix_i\, :\,  x\in\Z^n,\ Ax=0,\,
\sigma_kx_k\geq 1,\, 0\leq \sigma_ix_i\leq \sigma_ih_i,\, i=1,\dots,n\right\} .$$
The solution $x$ attaining minimum $\|x\|_1$ among these $n$ programs must be a
Graver element as desired. While we solve $n$ programs instead of one at each application
of Lemma \ref{sub_Graver}, all these programs are in standard form with matrix $A$. Since
$\Delta(A)\leq a$ which is fixed, the oracle algorithm requires the solution of polynomially
many programs over the same matrix $A$. Since $\Delta(A)\leq a$ implies in particular
that $\max|A_{i,j}|\leq a$, and it is assumed that $\min\{\td(A),\td(A^T)\}\leq d$,
Proposition \ref{sparse_IP} implies that each of these programs can be solved in polynomial time.
\epr

\section{Open problems}

It remains an important open problem whether a circuit walk $x^0,x^1,\dots,x^s=x^*$
could be computed in polynomial time {\em without} solving any linear program.
This might lead to a new algorithm for linear programming, and might shed light on the
fundamental long standing open problem of whether a {\em strongly polynomial time}
algorithm for linear programming exists, namely one performing a number of arithmetic
operations polynomial in $n$ and independent of the bit size of the data $A,w,b,l,u$.

We have shown in Theorem \ref{integer} that using an integer programming oracle we can
compute a Graver walk in time which depends polynomially on $\Delta(A)$ given in {\em unary}.
An interesting open problem is whether this could be strengthened to show that,
in fact, using such an oracle, a Graver walk could be computed in time which depends
polynomially on $\Delta(A)$ given in {\em binary}, that is, on the bit size of all data.

Another interesting open problem is whether the assumption on bounded tree-depth in
Theorem \ref{tree-depth} could be dropped, namely, is it true that for any fixed $a$,
Graver walks are polynomial time computable over any matrix $A$ with $\Delta(A)\leq a$.
A positive answer will imply a positive answer to the fundamental open problem
of whether integer programming could be solved in polynomial time over any matrix $A$
with $\Delta(A)\leq a$ (this is true for $a=1$ as then $A$ must be totally unimodular).

\section*{Acknowledgments}

The author is grateful to Jes\'{u}s A. De Loera, Sean Kafer and Laura Sanit\`{a} for useful
conversations on this work. He thanks the referees for useful suggestions which improved
the presentation. The research was partially supported by the National Science
Foundation under Grant No. DMS-1929284 while the author was in residence at the Institute
for Computational and Experimental Research in Mathematics, Brown University, during the
semester program {\em Discrete Optimization: Mathematics, Algorithms, and Computation},
and by the Dresner chair at the Technion.


\begin{thebibliography}{}

\bibitem{BDF}
Steffen Borgwardt, Jes\'{u}s A. De Loera, Elisabeth Finhold:
Edges versus circuits: a hierarchy of diameters in polyhedra.
Advacnces in Geometry 16:511--530, 2016.

\bibitem{BFH}
Steffen Borgwardt, Elisabeth Finhold, Raymond Hemmecke:
On the circuit diameter of dual transportation polyhedra.
SIAM Journal on Discrete Mathematics 29:113--121, 2015.

\bibitem{BV}
Steffen Borgwardt, Charles Viss:
A polyhedral model for enumeration and optimization over the set of circuits.
Discrete Applied Mathematics 308:68--83, 2022.

\bibitem{DHL}
Jes\'{u}s A. De Loera, Raymond Hemmecke, Jon Lee:
On augmentation algorithms for linear and integer-linear programming:
from Edmonds-Karp to Bland and beyond.
SIAM Journal on Optimization 25:2494--2511, 2015.

\bibitem{DKS}
Jes\'{u}s A. De Loera, Sean Kafer, Laura Sanit\`{a}:
Pivot rules for circuit-augmentation algorithms in linear optimization.
SIAM Journal on Optimization 32:2156--2179, 2022.

\bibitem{EHKKLO}
Friedrich Eisenbrand, Christoph Hunkenschroder, Kim-Manuel Klein,
Martin Kouteck\'y, Asaf Levin, Shmuel Onn:
Sparse integer programming is fixed-parameter tractable.
Mathematics of Operations Research, 16 pages, to appear,
{\tt https://doi.org/10.1287/moor.2023.0162}

\bibitem{HOW}
Raymond Hemmecke, Shmuel Onn, Robert Weismantel:
A polynomial oracle-time algorithm for convex integer minimization.
Mathematical Programming 126:97--117, 2011.

\bibitem{KLO}
Martin Kouteck\'y, Asaf Levin, Shmuel Onn:
A parameterized strongly polynomial algorithm for block structured integer programs.
Proceedings of ICALP 2018 (International Colloquium on Automata, Languages, and Programming),
Leibniz International Proceedings in Informatics, 107-85:1--14, 2018.

\bibitem{KO}
Martin Kouteck\'y, Shmuel Onn:
Sparse integer programming is FPT.
Bulletin of the European Association for Theoretical Computer Science 134:69--71, 2021.

\bibitem{NO}
Ne\v{s}et\v{r}il, J., Ossona de Mendez, P.:
Sparsity: Graphs, Structures, and Algorithms.
Algorithms and Combinatorics, Springer (2012)

\bibitem{Onn}
Shmuel Onn:
Nonlinear Discrete Optimization.
Zurich Lectures in Advanced Mathematics, European Mathematical Society (2010).
Available online,
{\tt https://sites.google.com/view/shmuel-onn/book}
\end{thebibliography}
\end{document}